\documentclass[10pt,twocolumn,a4paper]{article}

\usepackage[left=20mm, right=15mm, top=15mm, bottom=15mm]{geometry}

\usepackage{subfig}
\usepackage{flushend}
\usepackage{indentfirst}
\usepackage{graphics}
\usepackage{amsmath}
\usepackage{graphicx}
\usepackage{epstopdf}
\usepackage{float}

\usepackage[font=footnotesize,labelfont=bf]{caption}
\usepackage[labelsep=period]{caption}

\graphicspath{{./figure/}}

\newtheorem{ozn}{Definition}
\newtheorem{thm}{Theorem}
\newtheorem{zau}{Remark}
\newtheorem{lem}{Lema}
\newtheorem{nas}{Conclusion}

\begin{document}

\title{\huge \textbf{Robust extrapolation  problem for random processes with
stationary increments}}

\date{}

\twocolumn[

\begin{@twocolumnfalse}
\maketitle

Mathematics and Statistics 2(2): 78-88, 2014 \\
DOI: 10.13189/ms.2014.020204

\vspace{20pt}

\author{\textbf{Maksym Luz}, \textbf{Mikhail Moklyachuk}$^{*}$, \\\\
 {Department of Probability Theory, Statistics and Actuarial
Mathematics, \\
Taras Shevchenko National University of Kyiv, Kyiv 01601, Ukraine}\\
$^{*}$Corresponding Author: Moklyachuk@gmail.com}\\\\\\

\end{@twocolumnfalse}

]

\noindent \textbf{\Large{Abstract}} \hspace{2pt} The  problem  of optimal estimation of linear
functionals  $A {\xi}=\int_{0}^{\infty} a(t)\xi(t)dt$ and $A_T{\xi}=\int_{0}^{T} a(t)\xi(t)dt$ depending
on the unknown values of random process $\xi(t)$, $t\in R$, with
stationary $n$th increments from observations of ttis process for $t<0$ is considered.
Formulas for calculating mean square error and spectral
characteristic of optimal linear estimation of the functionals are
proposed in the case when spectral density is exactly known.
Formulas that determine the least favorable spectral densities are proposed for  given sets
of admissible spectral densities.\\

\noindent \textbf{\Large{Keywords}} \hspace{2pt} Random process with stationary increments;
minimax-robust estimate; mean square error; least favorable spectral
density; minimax-robust spectral characteristic\\

\noindent\hrulefill

\section{\Large{Introduction}}

Estimation of unknown values of random processes is an important part of the theory of random processes. A lot of researches were dedicated to the stationary case. Traditional methods of solution of the linear extrapolation,
interpolation and filtering problems for stationary stochastic
processes were developed by  Kolmogorov $\cite{Kolmogorov}$, Wiener$\cite{Wiener}$, Yaglom $\cite{Yaglom:1987a,Yaglom:1987b}$. The further results one can find in  book by Rozanov $\cite{Rozanov}$. Yaglom $\cite{Yaglom:1955,Yaglom:1957}$ generalized the case of stationary processes. He developed a theory of non-stationary processes whose increments of order $\mu\neq0$ define a stationary process. The spectral representation for stationary increments and canonical factorization for spectral densities were received, the problem of linear extrapolation of unknown value of stationary random increment from observation of the process was solved. Further results for such stochastic processes were
presented by Pinsker $\cite{Pinsk:1955}$, Yaglom and Pinsker $\cite{Pinsk:1954}$. See books by Yaglom $\cite{Yaglom:1987a,Yaglom:1987b}$ for more relative results and references.

The mean square optimal estimation problems for stochastic processes
with $n$th stationary increments are natural generalization of the
linear extrapolation, interpolation and filtering problems for
stationary stochastic processes.

 The classical methods of extrapolation, interpolation and filtering problems are based on the assumption
that the spectral density of the process is known.
In practice, however, it is impossible to obtain complete information
on the spectral density in most cases. To solve the problem one
finds parametric or nonparametric estimates of the unknown spectral
density or selects a density by other reasoning. Then
 the classical estimation method is applied provided that the estimated
or selected density is the true one.  Vastola and Poor $\cite{Vastola_Poor}$ have demonstrated that described procedure can result in
significant increasing of the value of error. This is a
reason for searching estimates which are optimal for all densities from
a certain class of the admissible spectral densities. These
estimates are called minimax since they minimize the maximal value
of the error. A survey of results in minimax (robust) methods of
data processing can be found in the paper by Kassam and Poor $\cite{Kassam_Poor}$.
The paper by  Grenander $\cite{Grenander}$ should be marked as the first one
where the minimax extrapolation problem for stationary processes was
formulated and solved. Franke and Poor$\cite{Franke}$, Franke $\cite{Franke:1985}$
investigated the minimax extrapolation and filtering problems for
stationary sequences with the help of convex optimization methods.
This approach makes it possible to find equations that determine the
least favorable spectral densities for various classes of admissible
densities. In papers
by Moklyachuk \cite{Moklyachuk:1989} - \cite{Moklyachuk:2001} the minimax approach was applied to
extrapolation, interpolation and filtering problems for functionals
which depend on the unknown values of stationary processes and
sequences. For more details see, for example, books by Kurkin et al.$\cite{Kurkin}$, Moklyachuk$\cite{Moklyachuk:2008}$. The case of vector stationary sequences and processes was developed by  Moklyachuk and Masyutka $\cite{Moklyachuk:2012}$. Dubovets'ka and Moklyachuk investigated periodically correlated stochastic sequences and random processes. In the articles $\cite{Dubovetska1,Dubovetska6}$ they considered the
minimax interpolation problem for the linear functionals which depend on unknown values of those sequences and processes.
Luz and Moklyachuk$\cite{Luz2012a,Luz2012b}$ solved the
minimax interpolation problem for the linear functional
$A_N\xi=\sum_{k=0}^Na(k)\xi(k)$ which depends on unknown (missed) values of a
stochastic sequence $\xi(m)$ with stationary $n$th increments from observations of the sequence with and without noise.

This article is dedicated to the mean square optimal estimates of the linear
functionals
\[A {\xi}=\int_{0}^{\infty}a(t)\xi(t)dt,\quad
A_T{\xi}=\int_{0}^Ta(t)\xi(t)dt\]
 which depend on the unknown values of a random process $\xi(t)$ with stationary $n$th increments. Estimates are
based on observations of the process $\xi(t)$ for
$t<0$. The estimation problem for processes with
stationary increments is solved in the case of spectral certainty
where the spectral density of the sequence is known as well as in
the case of spectral uncertainty where the spectral density of the
sequence is not known, but a set of admissible spectral densities is
given. Formulas are derived for computing the value of the
mean-square error and the spectral characteristic of the optimal
linear estimates of functionals $A{\xi}$ and $A_T{\xi}$ in the case
of spectral certainty where spectral density of the process is
known. Formulas that determine the least favorable spectral
densities and the minimax (robust) spectral characteristic of the
optimal linear estimates of the functionals are proposed in the case
of spectral uncertainty for concrete classes of admissible spectral
densities.

\section{\Large{Stationary random increment process. Spectral representation}}

\begin{ozn}
For a given random process $\{\xi(t),t\in  R\}$ a process
\begin{equation}
\label{oznachPryrostu_cont}
\xi^{(n)}(t,\tau)=(1-B_{\tau})^n\xi(t)=\sum_{l=0}^n(-1)^lC_n^l\xi(t-l\tau),
\end{equation}
where $B_{\tau}$ is a backward shift operator with step $\tau\in   R$, such that
$B_{\tau}\xi(t)=\xi(t-\tau)$, is called the random $n$th increment with step  $\tau\in  R$.
\end{ozn}

For the random $n$th increment process $\xi^{(n)}(t,\tau)$ the
following relations hold true:
\begin{equation}
\xi^{(n)}(t,-\tau)=(-1)^n\xi^{(n)}(t+n\tau,\tau), \label{tot1_cont}
\end{equation}
\begin{equation}
\label{tot2_cont}
\xi^{(n)}(t,k\tau)=\sum\nolimits_{l=0}^{(k-1)n}A_l\xi^{(n)}(t-l\tau,\tau),\quad
\forall k\in  N,
\end{equation}
where coefficients $\{A_l,l=0,1,2,\ldots,(k-1)n\}$ are
determined by the representation
\[(1+x+\ldots+x^{k-1})^n=\sum_{l=0}^{(k-1)n
}A_lx^l.
\]

\begin{ozn}
The random $n$th increment process $\xi^{(n)}(t,\tau)$ generated
by random process $\{\xi(t),t\in  R\}$ is wide sense
stationary if the mathematical expectations
\[
\mathsf E\xi^{(n)}(t_0,\tau)=c^{(n)}(\tau),
\]
\[
{\mathsf E}
\xi^{(n)}(t_0+t,\tau_1)\xi^{(n)}(t_0,\tau_2)=D^{(n)}(t,\tau_1,\tau_2)
\]
exist for all $t_0,\tau,t,\tau_1,\tau_2$ and do not depend on $t_0$. The function $c^{(n)}(\tau)$ is called the mean value of the $n$th increment and the function
$D^{(n)}(t,\tau_1,\tau_2)$ is called the structural function of the stationary $n$th increment (or the structural function of $n$th
order of the random process $\{\xi(t),t\in  R\}$).

The random process $\{\xi(t),t\in  R\}$ which determines the
stationary $n$th increment process $\xi^{(n)}(t,\tau)$ by formula
(\ref{oznachPryrostu_cont})  is called the process with stationary $n$th
increments.
\end{ozn}

\begin{thm}  \label{thm1_cont}
The mean value $c^{(n)}(\tau)$ and the structural function
$D^{(n)}(m,\tau_1,\tau_2)$ of a random stationary $n$th increment
process  $\xi^{(n)}(t,\tau)$ can be represented in the following
forms
\begin{equation}
\label{serFnaR_cont}
c^{(n)}(\tau)=c\tau^n,
\end{equation}
 \[D^{(n)}(t,\tau_1,\tau_2)=\]\begin{equation}\label{strFnaR_cont}=\int_{-\infty}^{\infty} e^{i\lambda
  t } (1-e^{-i\tau_1\lambda})^n(1-e^{i\tau_2\lambda})^n\frac{(1+ \lambda^2)^{n}}
{\lambda^{2n}}dF(\lambda),
\end{equation}
where  $c$ is a constant, $F(\lambda)$ is a left-continuous
nondecreasing bounded function with $F(-\infty)=0$. The constant $c$
and the function $F(\lambda)$ are determined uniquely by the
increment process $\xi^{(n)}(t,\tau)$.

From the other hand, a function $c^{(n)}(\tau)$ which has the form
$(\ref{serFnaR_cont})$ with a constant $c$ and a function
$D^{(n)}(m,\tau_1,\tau_2)$ which has the form $(\ref{strFnaR_cont})$ with a
function $F(\lambda)$ which satisfies the indicated conditions are
the mean value and the structural function of some stationary $n$th
increment process $\xi^{(n)}(t,\tau)$.
\end{thm}

Using representation $(\ref{strFnaR_cont})$ of the structural function
 of the stationary
$n$th increment process $\xi^{(n)}(t,\tau)$ and the Karhunen theorem
(see Karhunen $\cite{Karhunen}$), we get the following spectral representation of
the stationary $n$th increment process
$\xi^{(n)}(t,\tau)$:
\begin{equation}
\label{predZnaR_cont} \xi^{(n)}(t,\tau)=\int_{-\infty}^{\infty}
e^{it \lambda }(1-e^{-i\lambda\tau})^n\frac{(1+i\lambda )^n}{(i\lambda)^n}dZ(\lambda),
\end{equation}
where $Z(\lambda)$ -- is an orthogonal random measure on $R$ connected with the \emph{spectral function} $F(\lambda)$ by the relation
\begin{equation}
\label{FtaZ_cont}
 \mathsf EZ(A_1)\overline{Z(A_2)}=F(A_1\cap A_2)<\infty.
\end{equation}

Denote by $H(\xi^{(n)} )$ the subspace of the Hilbert space $H=L_2(\Omega,\mathcal{F}, \mathsf P)$ of the second order random variables which is generated by elements
$\{\xi^{(n)}( t ,\tau):t, \tau\in  R\}$ and let $H^{t}(\xi^{(n)} )$, $t \in   R$, be a  subspace of $H(\xi^{(n)} )$ generated by
elements
$\{\xi^{(n)}(u,\tau):u\leq t, \tau>0\}$.  Let $S { (\xi^{(n)})}$  is defined by relationship
\[S { (\xi^{(n)})}=\bigcap_{t\in  R} H^{t }(\xi^{(n)} ).\]
 Since the space $ S { (\xi^{(n)})}$ is a subspace in the Hilbert space $H (\xi^{(n)})$,
 the   space $H (\xi^{(n)})$ admits  the decomposition
  \[ H (\xi^{(n)} )= S { (\xi^{(n)})}\oplus R { (\xi^{(n)})},\]
  where $R { (\xi^{(n)})}$ is an orthogonal complement of the subspace
 $S { (\xi^{(n)})}$ in the space $H { (\xi^{(n)})}$.

From now we will consider increments $\xi^{(n)}(t,\tau)$ with step  $\tau>0$.

\begin{ozn} A stationary
$n$th increment process $\xi^{(n)}(t,\tau)$ is called regular if
 $H (\xi^{(n)} )=R { (\xi^{(n)})}$. It is called singular if
 $ H (\xi^{(n)} )=S { (\xi^{(n)})}$.
 \end{ozn}

\begin{thm}\label{thm 3_cont} A wide-sense stationary random increment process $\xi^{(n)}(t,\tau)$ admits a unique representation in the form
\begin{equation}\label{RozklVold_cont}\xi^{(n)}(t,\tau)=\xi_r^{(n)}(t,\tau)+\xi_s^{(n)}(t,\tau),\end{equation}
where $\{\xi_r^{(n)}(t,\tau):t\in  R\}$ is a regular increment
process and
$\{\xi_s^{(n)}(t,\tau):t\in  R\}$ is a singular
increment process. Moreover,  the increment process
$\xi_r^{(n)}(z,\tau)$ and $\xi_s^{(n)}(t,\tau)$ are orthogonal for all $t,z\in  R$.\end{thm}

Components of the representation $(\ref{RozklVold_cont})$ are constructed
in the following way: 
\[\xi_s^{(n)}(t,\tau)=\mathsf E[\xi^{(n)}(t,\tau)|S(\xi^{(n)})|,\]
\[
\xi_r^{(n)}(t,\tau)=\xi^{(n)}(t,\tau)-\xi_s^{(n)}(t,\tau).\]

Let $\{\eta(t):t\in  R\}$ be a random
process with independent increments such that $\mathsf E|\eta(t)-\eta(s)|^2=|t-s|$ and for all $z\in   R$  $H^{z}(\xi^{(n)} )=H^{z}(\eta )$, where  the subspace $H^{z}(\eta )$ of the space $H$ is generated by values  $\{ \eta(u):u\leq z \}$ of the process $\eta(t)$. Defined random process is called \emph{an innovate  process}.

\begin{thm}\label{thm 4_cont}A random stationary increment process
$\xi^{(n)}(t,\tau)$ is regular if and only if there exists an
innovate process $\{\eta(t):t\in  R\}$ and a function $\{\varphi^{(n)}(t,\tau):t\geq0\}$,
$\int_{0}^{\infty} |\varphi^{(n)}(t,\tau)|^2dt<\infty$, such that
\begin{equation}\label{odnostRuhSer_cont}\xi^{(n)}(t,\tau)=
\int_{0}^{\infty}\varphi^{(n)}(u,\tau)d\eta_u(t-u),\end{equation}
\end{thm}

\begin{nas} Using theorems $ \ref{thm 3_cont}$ and $\ref{thm 4_cont}$  one can conclude that a wide-sense stationary random
increment process admits a unique
representation in the form
\begin{equation}\xi^{(n)}(t,\tau)=\xi_s^{(n)}(t,\tau)+
\int_{0}^{\infty}\varphi^{(n)}(u,\tau)d\eta_u(t-u),\label{odnostRuhSerZag_cont}\end{equation} where
$\int_{0}^{\infty} |\varphi^{(n)}(t,\tau)|^2dt<\infty$ and $ \eta(t) $, $t\in   R$, is an innovate process.\end{nas}

Let the stationary $n$th increment process  $\xi^{(n)}(t,\tau)$
admit  the canonical representation $(\ref{odnostRuhSer_cont})$.
 In this case the spectral function
$F(\lambda)$ of the stationary increment process $\xi^{(n)}(t,\tau)$
has spectral density $f(\lambda)$ which admits the canonical
factorization
\begin{equation}\label{SpectrRozclad_f_cont}f(\lambda)= |\Phi({\lambda})|^2,\quad\Phi (\lambda)=\int_{0}^{\infty} e^{-i\lambda t}\varphi(t)dt.
\end{equation} Let us define
$$\Phi_{\tau}(\lambda)=\int_{0}^{\infty} e^{-i\lambda t}\varphi^{(n)}(t,\tau) dt=\int_{0}^{\infty}e^{-i\lambda t}\varphi_{\tau}(t)dt,$$ where
$\varphi_{\tau}(t)= \varphi^{(n)}(t,\tau)$ is the function from the representation $(\ref{odnostRuhSer_cont})$. Defined function $\Phi_{\tau}(\lambda)$, which is a Fourier transform of the function $\varphi^{(n)}(t,\tau)$,  is related  with spectral density $f(\lambda)$  of the random process  $\xi^{(n)}(t,\tau)$ by relations
\begin{equation}\left|\Phi_{\tau}({\lambda})\right|^2=\frac{|1-e^{-i\lambda\tau}|^{2n}(1+\lambda^2)^n}{\lambda^{2n}}f(\lambda),
\label{dd_cont}\end{equation}
\begin{equation}\Phi_{\tau}({\lambda})=\frac{(1-e^{-i\lambda\tau})^n(1+i\lambda)^n}{(i\lambda)^{n}}\Phi(\lambda).
\label{ddd_cont}\end{equation}

 The one-sided moving average representation (\ref{odnostRuhSer_cont}) is used for finding the optimal mean square
estimate of the unknown values of a process $\xi(t)$ from known observation for $t<0$.

\section{\Large{Hilbert space projection method of extrapolation of linear functionals}}

Let a random process $\{\xi(t),t\in  R\}$ defines $n$th increment $\xi^{(n)}(t,\tau)$ with an absolutely
continuous spectral function $F(\lambda)$
 which has spectral density $f(\lambda)$. Without loss of generality we will assume that the mean value of the increment process $\xi^{(n)}(t,\tau)$ equals to 0.
 Let the stationary increment process $\xi^{(n)}(t,\tau)$
admit the one-sided moving average representation
 $(\ref{odnostRuhSer_cont})$ and the spectral density  $f(\lambda)$ admits
the canonical factorization $(\ref{SpectrRozclad_f_cont})$. Consider the
case where the step $\mu>0$. Let we know the values of the process $\xi(t)$ for $t<0$. The problem is to find the mean square optimal linear
estimates of  functionals
$A_T\xi=\int_{0}^{T} a(t)\xi(t)dt$ and $A\xi=\int_{0}^{\infty} a(t)\xi(t)dt$ which depend on unknown values
 $\xi(t)$, $t\geq0$.

 In order to solve the stated problem we will present the process   $\xi(t)$, $t\geq0$, as a sum of some its increments  $\xi^{(n)}(t,\tau)$,  $t\geq0$, $\tau>0$, and some of its initial values $\xi^0$. Particularly, when  $\tau^*>t^*\geq0$, a relation   \[\xi(t^*)=\xi^{(n)}(t^*,\tau^*)+\sum_{l=1}^n(-1)^lC_n^l\xi(t^*-l\tau^*)\] comes from $(\ref{oznachPryrostu_cont})$, where $\xi^0=\{\xi(t^*-l\tau^*):l=1,2,\ldots,n\}\subset\{\xi(t):t\leq0\}$ are known observations. The following lema describes a representation of the functional $A\xi$ from some of known initial values of the process $\xi(t)$ and its increments $\xi^{(n)}(t,\tau)$ for $tgeq0$ in  case of arbitrary step $\tau>0$.

\begin{lem}\label{lema predst A_cont}
A linear functional $A\xi=\int_{0}^{\infty}a(t)\xi(t)dt$ admits a representation  $A\xi=B \xi-V\xi$, where
\[B\xi=\int_{0}^{\infty} b_{\tau}(t)\xi^{(n)}(t,\tau)dt,\quad
V\xi=\int_{-\tau n}^{0}v_{\tau}(t)\xi(t)dt,
\]
\begin{equation} \label{koefv_cont}
v_{\tau}(t)=\sum_{l=\left[-\frac{t}{\tau}\right]'}^n(-1)^lC_n^lb_{\tau}(t+l\tau),\quad
t\in [-\tau n;0), \end{equation}
\begin{equation} \label{koef b_cont}b_{\tau}(t)=\sum_{k=0}^{\infty}a( t +\tau k)d(k)=\mathbf D^{\tau}a(t) ,\,\,t\geq0,
\end{equation}
where $[x]'$ denotes the least integer number among
numbers greater or equal to $x$, $\{d(k):k\geq0\}$ are coefficients determined
by the relation
$\sum_{k=0}^{\infty}d (k)x^k=\left(\sum_{j=0}^{\infty}x^{
j}\right)^n$, $\mathbf D^{\tau}$ is a linear transformation which acts on an arbitrary function $x(t)$, $t>0$, by formula \begin{equation}\mathbf
D_{\tau}x(t)=\sum_{k=0}^{\infty}x(t+\tau k )d(k).\label{operator D_cont}\end{equation}
\end{lem}

\emph{Proof.} From $(\ref{oznachPryrostu_cont})$ we can obtain the formal equation
\begin{equation}
\label{obernene oznach prir_cont}
\xi(t)=\frac{1}{(1-B_{\tau})^n}\xi^{(n)}(t,\tau)=
\sum_{j=0}^{\infty}d(j)\xi^{(n)}( t- \tau j,\tau),
\end{equation}
which follows the relations
\[
\int_{0}^{\infty} a(t)\xi(t)dt=-\int_{-\tau n}^0v_{\tau}(t) \xi(t)dt+\]\[+
\int_{0}^{\infty} \left(\sum_{k=0}^{\infty}a(t+\tau k )d(k)\right)\xi^{(n)}(t,\tau)dt,
\]
\[
\int_{0}^{\infty} b_{\tau}(k)\xi^{(n)}(t,\tau)dt=\]\[= \int_{ -\tau
n}^{0}\xi(t)\sum_{l=\left[-\frac{t}{\tau}\right]'}^{n}(-1)^lC_n^k b_{\tau}(t+\tau l)dt+\]
\[+
\int_{ 0}^{\infty}\xi(t) \sum_{l=0}^{n}(-1)^lC_n^k b_{\tau}(t+\tau l) dt.
\]
From  two of last relations we can get the representation of the functional
 $A\xi$ and relations $(\ref{koef b_cont})$,
$(\ref{koefv_cont})$.

\begin{nas}
The linear functional $A_T\xi$ admits a representation $A_T\xi=B_T \xi-V_T\xi$, where
\[B_T\xi=\int_0^T b_{\tau,T}(t)\xi^{(n)}(t,\tau)dt,\,
V_T\xi=\int_{-\tau n}^{0}v_{\tau,T}(t)\xi(t)dt,
\] and functions $b_{\tau,T}(t)$, $t\in[0,T]$, and  $v_{\tau,T}(t)$, $t\in[-\tau n;0)$, are defined by formulas  $(\ref{koefv_cont})$ and
$(\ref{koef b_cont})$ respectively stating $a(t)=0$ when $t>T$.
\end{nas}

We will suppose that the following conditions on the function $b_{\tau}(t)$ hold true
\begin{equation}\label{umova na b ex_cont}\int_{0}^{\infty}|b_{\tau}(t)|dt<\infty,\quad
\int_{0}^{\infty} t|b_{\tau}(t)|^2dt<\infty.\end{equation}
 Under these conditions the functional $B\xi$ has the second moment and
the operator $\mathbf B^{\tau}$  defined below is
 compact.
Since the functions $a(t)$ and $b_{\tau}(t)$ are related by  $(\ref{koef b_cont})$,
 the following conditions hold true
\begin{equation}\label{umova na a
ex_cont}\int_0^{\infty}|\mathbf D^{\tau}a(t)|dt<\infty,\quad
\int_0^{\infty} t|\mathbf D^{\tau}a(t)|^2dt<\infty.\end{equation}

Let $\widehat{A}\xi$ denote the mean square optimal linear estimate
of the functional $A\xi$ from observations of the process $\xi(t)$
for
 $t<0$ and let
$\widehat{B}\xi$ denote the mean square optimal linear estimate of
the functional
 $B\xi$ from observations of the random $n$th increment process $\xi^{(n)}(t,\tau)$ for
$ t<0$. Let $\Delta(f,\widehat{A}\xi):=\mathsf{E}
|A\xi-\widehat{A}\xi|^2$ denote the mean square error of the
estimate $\widehat{A}\xi$ and let
$\Delta(f,\widehat{B}\xi):=\mathsf{E} |B\xi-\widehat{B}\xi|^2$
denote the mean square error of the estimate $\widehat{B}\xi$. Since
values
$\xi(t)$ for $t\in [-\tau n;0)$ are known, the following
equality comes from lema $\ref{lema predst A_cont}$:
\begin{equation}\label{main
formula_cont}\widehat{A}\xi=\widehat{B}\xi-V\xi.\end{equation}
 Thus
\[\Delta(f,\widehat{A}\xi)=\mathsf E |A\xi-\widehat{A}\xi|^2=\mathsf E|A\xi+V\xi-\widehat{B}\xi|^2=\]\[=\mathsf E|B\xi-\widehat{B}\xi|^2=\Delta(f,\widehat{B}\xi).\]

Denote by $L_2^{0-}(f)$ the subspace of the Hilbert space $L_2(f)$
generated by the set of functions
 \[h(\lambda)=(1-e^{-i\lambda\tau})^{n}\dfrac{(1+i\lambda)^{n}}{(i\lambda)^{n}}\int_0^{\infty} h(t)e^{-i\lambda t}dt.\]
Every
linear estimate $\widehat{B}\xi$ of the functional $B\xi$ admits a representation
\begin{equation} \label{est1_cont} \widehat{B}\xi=\int_{-\infty}^{\infty}
h_{\tau}(\lambda)(1-e^{-i\lambda\tau})^{n}\dfrac{(1+i\lambda)^{n}}{(i\lambda)^{n}}dZ(\lambda), \end{equation}
 where $h_{\tau}(\lambda)$ is the spectral characteristic of the estimate $\widehat{B}\xi$. The spectral characteristic of the optimal estimate
provides the minimum value of the mean square error
$\Delta(f,\widehat{B}\xi)$.

  Let the random increment $\xi^{(n)}(t,\tau)$ admits the canonical representation
$(\ref{odnostRuhSer_cont})$. Then the functional $B\xi$ can be presented as
\[B\xi=\int_0^{\infty}\int_0^{\infty}b_{\tau}(t)\varphi(u,\tau)d\eta_u(t-u)dt=\]
\[=\int_{-\infty}^0\int_0^{\infty}b_{\tau}(t)\varphi(t-u,\tau)dt d\eta(u)+\]\[+\int_0^{\infty}\int_{u}^{\infty}b_{\tau}(t)\varphi(t-u,\tau)dtd\eta(u).\] As the relation  $H^{0}(\xi^{(n)} )=H^{0}(\eta )$ holds true and increments of the process $\eta(t)$ are orthogonal, the optimal estimate $\widehat{B}\xi$ of the functional $B\xi$ is calculated as
\[\widehat{B}\xi=\int_{-\infty}^0\int_0^{\infty}b_{\tau}(t)\varphi_{\tau}(t-u )dtd\eta(u)=\]\begin{equation}\label{pobudova otsinky_1_cont}+\int_{-\infty}^0B\varphi_{\tau}(t-u )d\eta(u)=\int_{-\infty}^{\infty}B^*\varphi_{\tau}(\lambda )d\eta^*(\lambda),\end{equation}
where $B^*\varphi_{\tau}(\lambda )$ and $\eta^*(\lambda)$ are inverse Fourier transforms of the function
$B\varphi_{\tau}(u )=\int_0^{\infty }b_{\tau}(t)\varphi_{\tau}(t-u )dt$, $u<0$, and the process $\eta(u)$ respectively.
Yaglom\cite{Yaglom:1955} showed that
 \begin{equation}\label{pobudova otsinky_2_cont}\eta ^*(\lambda)=\int_{-\infty}^{\lambda}\frac{dZ(p)}{\Phi(p)}.\end{equation} So we need to find $B^*\varphi_{\tau}(\lambda )$.
\[B^*\varphi_{\tau}(\lambda )=\int_{-\infty}^0e^{i\lambda s}\int_0^{\infty}b_{\tau}(t)\varphi_{\tau}(t-s )dt ds=\]\[=\int_0^{\infty}e^{i\lambda t}b_{\tau}(t)\int_0^{\infty}e^{-i\lambda(s+t)}\varphi_{\tau}(t+s )dsdt=\]\[=\int_0^{\infty}e^{i\lambda t}b_{\tau}(t)\int_t^{\infty}e^{-i\lambda z}\varphi_{\tau}(z )dzdt=\]\[=B_{\tau}(\lambda)\Phi_{\tau}({\lambda})-\int_0^{\infty}e^{i\lambda t}b_{\tau}(t)\int_0^{t}e^{-i\lambda z}\varphi_{\tau}(z )dzdt=\]\begin{equation}\label{pobudova otsinky_3_cont}=B_{\tau}(\lambda)\Phi_{\tau}({\lambda})-\int_0^{\infty}e^{i\lambda y}\int_0^{\infty}b_{\tau}(y+z)\varphi_{\tau}(z)dzdy.\end{equation}
Substituting the expressions $(\ref{pobudova otsinky_2_cont})$ and $(\ref{pobudova otsinky_3_cont})$ in $(\ref{pobudova otsinky_1_cont})$ and using $(\ref{ddd_cont})$ one can obtain the following formulas for culculating the spectral characteristic of the optimal estimate $\widehat{B}\xi$:
\begin{equation}\label{spectr_har_B_cont}h_{\tau}(\lambda)=
B^\tau( \lambda )-r_{\tau}( \lambda )
\Phi_{\tau}^{-1}({\lambda}),\end{equation}
\[B^\tau( \lambda )=\int_{0}^{\infty} b_\tau(t) e^{i\lambda t}dt,\,  r_{\tau}( \lambda )=
\int_{0}^{\infty} e^{i\lambda t} (\mathbf B^\tau \varphi_{\tau})( t)dt ,\] where
$\mathbf B^{\tau}$ is a linear operator in  $L_2([0,\infty))$ space which defined as
\[(\mathbf B^{\tau} \varphi_{\tau})(t)=\int_{0}^{\infty} b_\tau(t+u)\varphi_{\tau} (u )du.\]  Here $\varphi_{\tau} (u )=\varphi^{(n)}(u,\tau)$ is the function from moving average representation $(\ref{odnostRuhSer_cont})$. The operator $\mathbf B^{\tau}$ is compact providing  $(\ref{umova na b ex_cont})$.

The value of the mean square error
$\Delta(f,\widehat{B}\xi)$ can be calculated by the formula
\[\Delta(f,\widehat{B}\xi)=\mathsf E|B \xi-\widehat{B} \xi|^2=\]\begin{equation}\label{pohybka B_cont}=
\frac{1}{2\pi}\int_{-\infty}^{\infty}|r_{\tau}( \lambda )|^2
d\lambda=
||\mathbf B^{\tau}\varphi_{\tau}||^2.\end{equation}

\begin{thm} \label{thm2_cont}
Let a random process $\{\xi(t),m\in    Z\}$ determine a
stationary random $n$th increment process $\xi^{(n)}(t,\tau)$
with absolutely continuous spectral function $F(\lambda)$ and
spectral density $f(\lambda)$
 which admits the canonical factorization $(\ref{SpectrRozclad_f_cont})$.
The optimal linear estimate $\widehat{B}\xi$ of the functional
$B\xi$ which depends on the unobserved values $\{\xi^{(n)}(t,\tau):
t\geq0\}$, $\tau>0$, from observations of the process $\xi(t)$ for $t<0$  can be calculated by formula (\ref{est1_cont}).
The spectral characteristic $h_{\mu}(\lambda)$ of the
optimal linear estimate $\widehat{B}\xi$ can be calculated by
formula (\ref{spectr_har_B_cont}).
The value of the mean square error
$\Delta(f,\widehat{B}\xi)$ can be calculated by formula (\ref{pohybka B_cont}).
\end{thm}

Using  Theorem $\ref{thm2_cont}$ and representation $(\ref{RozklVold_cont})$, we can obtain an optimal estimate of an unobserved value $\xi^{(n)}(u,\tau)$, $\tau>0$, in the point $ u\geq0$ from observations of the
process $\xi(t)$ for $t<0$. The singular component $\xi_s^{(n)}(u,\tau)$ decomposition $(\ref{RozklVold_cont})$ of the process has errorless estimate. We will
use formula $(\ref{spectr_har_B_cont})$ to obtain the spectral
characteristic $h_{u,\tau}(\lambda)$ of the optimal estimate
$\widehat{\xi}^{(n)}(t,\tau)$ of the regular component
$\xi_r^{(n)}(u,\tau)$ of the process. Consider a function $B^\tau( \lambda )=e^{i\lambda u} $. It
follows from the derived formulas that the spectral characteristic
of the estimate
\[\widehat{\xi}^{(n)}(u,\tau)=\xi_s^{(n)}(u,\tau)+\]\begin{equation}\label{est2_cont}+\int_{-\infty}^{\infty}
h_{u,\tau}(\lambda)(1-e^{-i\lambda\tau})^n\dfrac{(1+i\lambda)^n}{(i\lambda)^n}dZ(\lambda)\end{equation}
 can be calculated by the formula
\begin{equation}\label{spectr_har_delta_xi_cont}h_{u,\tau}(\lambda)=e^{i\lambda u}-
 \Phi_{\tau}^{-1}({\lambda})\int_0^u\varphi_{\tau}(y)e^{-i\lambda
y}dy.\end{equation} The value of the mean square error can be calculated by the formula
\begin{equation}\label{pohybka delta_xi_cont}\Delta(f,\widehat{\xi}^{(n)}(u,\tau))=
 \frac{1}{2\pi}\int_{0}^{u }|\varphi_{\tau}(y)|^2dy.\end{equation}

The following statement holds true.

\begin{nas} \label{nas1_cont}
The optimal linear estimate $\widehat{\xi}^{(n)}(u,\tau)$ of the
value $\xi^{(n)}(u,\tau)$, $\tau>0$, in the point $u\geq0$ of the increment process $\xi^{(n)}(t,\tau)$
from observations of the process  $\xi(t)$, $t<0$, can be calculated by  formula (\ref{est2_cont}).
The
spectral characteristic  $h_{u,\tau}(\lambda)$ of the optimal linear
estimate $\widehat{\xi}^{(n)}(u,\tau)$ can be calculated by formula
$(\ref{spectr_har_delta_xi_cont})$. The value of mean square error
$\Delta(f,\widehat{\xi}^{(n)}(u,\tau))$ of the optimal linear
estimate can be calculated by formula (\ref{pohybka delta_xi_cont}).
\end{nas}

 Making use relation $(\ref{main
formula_cont})$ we can find the optimal
estimate $\widehat{A}\xi$ of the functional $A\xi$ from observations
of the process $\xi(t)$ for $t<0$. These estimate can be presented in the following form:
\[\widehat{A}\xi=-\int_{-\tau n}^{0}v_{\tau}(t)\xi(t)dt+\]\begin{equation} \label{estim A_cont}+\int_{-\infty}^{\infty}
h_{\tau}^{(a)}(\lambda)(1-e^{-i\lambda\tau})^n\dfrac{(1+i\lambda)^n}{(i\lambda)^n}dZ(\lambda),\end{equation}  where the function $v_{\tau}(t)$, $t=[-\tau n,0)$, is defined by relation $(\ref{koefv_cont})$. Using the relationship $(\ref{koef b_cont})$ between the functions $a(t)$ and $b_{\tau}(t)$, $t\geq0$, we obtain the following equation:
 \[(\mathbf B^{\tau}\mathbf \varphi_{\tau})(t)=\int_{0}^{\infty} \mathbf D^{\tau}a(t+u)\varphi (u,\tau)du=\]\[=\mathbf D^{\tau}\int_{0}^{\infty} a(t+u)\varphi (u,\tau)du=\mathbf D^{\tau} (\mathbf A\varphi_{\tau})(t),\] where the linear operator $\mathbf A$ is defined by the function $a(t)$, $t\geq0$, in the following way:
 \[(\mathbf A\varphi_{\tau})(t)=\int_{0}^{\infty} a(t+u)\varphi_{\tau} (u )du.\]
Thus the spectral characteristic and the value of the
mean square error of the optimal estimate $\widehat{A}\xi$ can be
calculated by the formulas
\begin{equation}\label{spectr_har_A_cont}h_{\tau}^{(a)}(\lambda)=A_{\tau}(\lambda)-r_{\tau}^{(a)}(\lambda)
\Phi_{\tau}^{-1}(\lambda),\end{equation}
\[A_{\tau}(\lambda)=\int_{0}^{\infty}\mathbf D^{\tau}a(t)
e^{i\lambda t}dt,\] \begin{equation}\label{r_A_cont}r_{\tau}^{(a)}(\lambda)=
\int_{0}^{\infty}\mathbf D^{\tau} (\mathbf A\varphi_{\tau})(t)e^{i\lambda t}dt.\end{equation}
\[\Delta(f,\widehat{A}\xi)=\mathsf E|A \xi-\widehat{A} \xi|^2=\]\begin{equation}\label{pohybka A_cont}=
\frac{1}{2\pi}\int_{-\infty}^{\infty}|r_{\tau}^{(a)}(\lambda)|^2
d\lambda=
||\mathbf D^{\tau}\mathbf A\varphi_{\tau}||^2.\end{equation}

The following theorem holds true.

\begin{thm} \label{thm3_cont}
Let a random process  $\{\xi(t),t\in  R\}$ determine
a stationary random $n$th increment process $\xi^{(n)}(t,\tau)$
with absolutely continuous spectral function $F(\lambda)$ and
spectral density $f(\lambda)$
 which admits the canonical factorization (\ref{SpectrRozclad_f_cont}).
The optimal linear estimate $\widehat{A}\xi$ of the functional
$A\xi$ of unobserved values $\xi(t)$,
$t\geq0$, from
observations of the process $\xi(t)$ for $t<0$ can be
calculated by formula (\ref{estim A_cont}). The spectral characteristic
$h_{\mu}^{(a)}(\lambda)$ of the optimal linear estimate
$\widehat{A}\xi$ can be calculated by formula (\ref{spectr_har_A_cont}).
The value of the mean square error $\Delta(f,\widehat{A}\xi)$ of the
optimal linear estimate can be calculated by formula (\ref{pohybka
A_cont}).
\end{thm}

Consider now the problem of the mean square optimal estimation of
the functional $A_T\xi$. The optimal estimate of the
functional can be calculated by formula
\[\widehat{A}_T\xi=-\int_{-\tau n}^{0}v_{\tau,T}(t)\xi(t)dt+\]\begin{equation} \label{estim A_T_cont}+\int_{-\infty}^{\infty}
h_{\tau,T}^{(a)}(\lambda)(1-e^{-i\lambda\tau})^n\dfrac{(1+i\lambda)^n}{(i\lambda)^n}dZ(\lambda),\end{equation}  where the function $v_{\tau,T}(t)$, $t\in[-\tau n;0)$, can be calculated by formulas
\[v_{\tau,T}(t)=\sum_{l=\left[-\frac{t}{\tau}\right]'}^{\min\left\{\left[\frac{T-t}{\tau}\right],n\right\}}(-1)^lC_n^lb_{\tau, T}(l\tau+t),\, t\in[-\tau n;0),\]\[ b_{\tau,T}(t)=\sum_{k=0}^{\left[\frac{T-t}{\tau}\right]}a( t +\tau k)d(k)=\mathbf D^{\tau}_Ta(t),\,
t\in[0;T].\] Here $\mathbf D^{\tau}_T$ is a linear transformation  which acts on an arbitrary function $x(t)$, $t\in[0,T]$, as \[\mathbf D^{\tau}_Tx(t)=\sum_{k=0}^{\left[\frac{T-t}{\tau}\right]}x( t +\tau k)d(k).\]

The spectral characteristic $h_{\tau,T}^{(a)}(\lambda)$ and the value of the mean square error  $\Delta(f,\widehat{A}_T\xi)$ of the estimate $\widehat{A}_T\xi$ can be calculated by  formulas
\begin{equation}\label{spectr_har_A_T_cont}h_{\tau,T}^{(a)}(\lambda)=A_{\tau, T}(\lambda)-r_{\tau,T}^{(a)}(\lambda)
\Phi_{\tau}^{-1}(\lambda),\end{equation}
\[A_{\tau, T}(\lambda)=\int_{0}^{T}\mathbf D^{\tau}_Ta(t)
e^{i\lambda t}dt,\] \begin{equation}\label{r_A_T_cont}r_{\tau,T}^{(a)}(\lambda)=
\int_{0}^{T}\mathbf D^{\tau}_T (\mathbf A_T\varphi_{\tau})(t)e^{i\lambda t}dt,\end{equation} where $\mathbf A_T $ is a linear operator in $L_2([0,\infty))$ space  defined by formula
\[(\mathbf A_T\varphi_{\tau})(t)=\int_{0}^{T-t} a(t+u)\varphi_{\tau}(u)du,\] and linear opperator $\mathbf D^{\tau}_T \mathbf A_T\varphi_{\tau}$ in $L_2([0,\infty))$ space is defined by formula \[\mathbf D^{\tau}_T (\mathbf A_T\varphi_{\tau})(t)=\sum_{k=0}^{\left[\frac{T-t}{\tau}\right]}\int_0^{T-t-\tau k}a(u+t+\tau k)\varphi_{\tau}(u)d(k)du;\]
\[\Delta(f,\widehat{A}_T\xi)=\mathsf E|A_T\xi-\widehat{A}_T \xi|^2=\]\begin{equation}\label{pohybka A_T_cont}=
\frac{1}{2\pi}\int_{-\infty}^{\infty}|r_{\tau,T}^{(a)}(\lambda)|^2
d\lambda=
||\mathbf D^{\tau}_T\mathbf A_T\varphi_{\tau}||^2.\end{equation}

Consequently, the following theorem holds true.

\begin{thm} \label{thm4_cont}
Let a random process $\{\xi(t), t\in   R\}$ determine
a stationary random $n$th increment process $\xi^{(n)}(t,\tau)$
with absolutely continuous spectral function $F(\lambda)$ and
spectral density $f(\lambda)$
 which admits the canonical factorization (\ref{SpectrRozclad_f_cont}).
The optimal linear estimate $\widehat{A}_T\xi$ of the functional
$A_T\xi$ of unobserved values $\xi(t)$, $t\geq0$, from
observations of the process $\xi(t)$ for $t<0$  can be
calculated by formula (\ref{estim A_T_cont}). The spectral characteristic
$h_{\tau,T}^{(a)}(\lambda)$ of the optimal linear estimate
$\widehat{A}_T\xi$ can be calculated by formula
(\ref{spectr_har_A_T_cont}). The value of mean square error
$\Delta(f,\widehat{A}_T\xi)$ can be calculated by formula
(\ref{pohybka A_T_cont}).
\end{thm}

Consider the case where $\tau>u\geq0$. In this case the optimal mean
square estimate of the value $\xi(u)$ in the point $u\geq0$ from observations $\xi(t)$ for $t<0$ can be calculated
by formula \[\widehat{\xi}(u)=\sum_{l=1}^T(-1)^{l+1}C_n^l\xi(u-l\tau)+\]\begin{equation}\label{estim_xi_cont}+\int_{-\infty}^{\infty}
h_{u,\tau}(\lambda)(1-e^{-i\lambda\tau})^n\dfrac{(1+i\lambda)^n}{(i\lambda)^n}dZ(\lambda).\end{equation}
The spectral characteristic $h_{u,\tau}(\lambda)$ and the
value of the mean square error $\Delta(f,\widehat{\xi}(u))=\Delta(f,\widehat{\xi}^{(n)}(u,\tau))$ of the estimate $\widehat{\xi}(u)$ can be calculated by
formulas $(\ref{spectr_har_delta_xi_cont})$ and $(\ref{pohybka delta_xi_cont})$ respectively.

Consequently, the following statement holds true.

\begin{nas} \label{nas2_cont} Let $\tau>u\geq0$. The optimal mean square estimate $\widehat{\xi}(u)$ of the element
$\xi(u)$ from observations of the process $\xi(t)$
for $t<0$ can be calculated by formula (\ref{estim_xi_cont}).
The spectral characteristic $h_{u,\tau}(\lambda)$ of the optimal
linear estimate $\widehat{\xi}(u)$ can be calculated by formula
(\ref{spectr_har_delta_xi_cont}). The value of mean square error
$\Delta(f,\widehat{\xi}(u))$ can be calculated by formula
(\ref{pohybka delta_xi_cont}).
\end{nas}

 \begin{zau}
Using relation $(\ref{dd_cont})$ we can find a relationship between
functions $\varphi_{\tau}(t)$, $t\geq0$,
 and $\varphi(t)$, $t\geq0$. So far as \[\int_{-\infty}^{\infty}\left|\ln \frac{|1-e^{-i\lambda\tau}|^{2n}(1+\lambda^2)^n}{\lambda^{2n}}\right|\frac{1}{1+\lambda^2}d\lambda <\infty\] for all $n\geq1$ and $\tau>0$, there are functions $\omega_{\tau}(t)$, $t\geq0$, and  \[\Omega_{\tau}(\lambda)=\int_{0}^{\infty}\omega_{\tau}(\lambda)e^{-i\lambda t}dt\] such that \[||\omega_{\tau}(\lambda)||^2=\frac{1}{2\pi}\int_{0}^{\infty}|\omega_{\tau}(\lambda)|^2dt<\infty,\] \[\frac{|1-e^{-i\lambda\tau}|^{2n}(1+\lambda^2)^n}{\lambda^{2n}}=|\Omega_{\tau}(\lambda)|^2\] and the following relation holds true:
 \begin{equation} \Phi_{\tau}(\lambda) =
 \Omega_{\tau}(\lambda)\Phi (\lambda).
 \label{dddd_cont}\end{equation}
 From $(\ref{dddd_cont})$ using inverse Fourier transform  we get  \[\varphi_{\tau}(t)= \int_{-\infty}^{\infty}e^{i\lambda t}\Omega_{\tau}(\lambda)\Phi(\lambda)d\lambda=\]\[=\int_0^{\infty}\varphi(x)\int_{-\infty}^{\infty}e^{i\lambda(t-x)}
 \Omega_{\tau}(\lambda)d\lambda dx=\]\[=\int_0^t\omega_{\tau}(t-x)\varphi(x)dx.\]
 Therefore, the functions $\varphi_{\tau}(t)$, $t\geq0$, and $\varphi(t)$, $t\geq0$, from $L_2([0,\infty))$ space are related by the relation
  \begin{equation}\varphi_{\tau}(t) = \mathbf W^{\tau}\varphi(t)=\int_0^t\omega_{\tau}(t-x)\varphi(x)dx,  \label{dd1_cont}\end{equation}
 where   $\mathbf W^{\tau}$ is a linear operator in $L_2([0,\infty))$ space defined by the function $\omega_{\tau}(t)$, $t\geq0$, from $L_2([0,\infty))$ space.
 When the functions $\varphi_{\tau}(t)$ and
 $\varphi(t)$ are defined on segment $[0,T]$, which means $\varphi_{\tau}(t)=\varphi (t)=0$ for $t>T$, the relation between them is defined by $(\ref{dd1_cont})$ for $t\in[0,T]$.
 \end{zau}

\section{\Large{Minimax-robust method of extrapolation}}
The proposed formulas may be employed under the condition that the
spectral density $f(\lambda)$ of the considered random process
$\xi(t)$ with stationary $n$th increments is known. The value of the
mean square error
$\Delta(h_{\tau}^{(a)}(f);f):=\Delta(f,\widehat{A}\xi)$ and the
spectral characteristic $h_{\tau}^{(a)}(f)$  of the optimal linear
estimate $\widehat{A}\xi$ of the functional $A\xi$ which depends of
unknown values  $\xi(t)$ can be calculated by  formulas
$(\ref{spectr_har_A_cont})$ and $(\ref{pohybka A_cont})$, the value of mean
square error
$\Delta(h_{\tau,T}^{(a)}(f);f):=\Delta(f,\widehat{A}_T\xi)$ and the
spectral characteristic $h_{\tau,T}^{(a)}(f)$ of the optimal linear
estimate $\widehat{A}_T\xi$
 of the functional $A\xi$ which depends of unknown values  $\xi(t)$, $t\geq0$,  can be calculated by  formulas
 $(\ref{spectr_har_A_T_cont})$ and $(\ref{pohybka A_T_cont})$.
In the case where the spectral density is not known, but a set $\mathcal
D$ of admissible spectral densities is given, the minimax (robust)
approach to estimation of the functionals of the unknown values of a
random process with stationary increments is reasonable. In
other words we are interesting in finding an estimate that minimizes
the maximum of the mean square errors for all spectral densities
from a given class $\mathcal   D$ of admissible spectral densities
simultaneously.

\begin{ozn} For a given class of spectral densities $\mathcal{D}$ a
spectral density $f_0(\lambda)\in\mathcal{D}$ is called least
favorable in $\mathcal{D}$ for the optimal linear estimate the
functional $A \xi$ if the following relation holds true:
\[\Delta(f^0)=\Delta(h_{\tau}^{(a)}(f^0);f^0)=\max_{f\in\mathcal{D}}\Delta(h_{\tau}^{(a)}(f);f).\]
\end{ozn}

\begin{ozn} For a given class of spectral densities $\mathcal{D}$
a spectral characteristic $h^0(\lambda)$ of the optimal linear estimate of the functional
$A \xi$ is called minimax-robust if there are satisfied conditions
\[h^0(\lambda)\in
H_{\mathcal{D}}=\bigcap_{f\in\mathcal{D}}L_2^{0-}(f), \]
\[\min_{h\in H_{\mathcal{D}}}\max_{f\in
\mathcal{D}}\Delta(h;f)=\sup_{f\in\mathcal{D}}\Delta(h^0;f).\]
\end{ozn}

Analyzing the derived formulas and using the introduced definitions
we can conclude that the following statements are true.

\begin{lem} Spectral density $f^0(\lambda)\in\mathcal{D}$ which admits the
canonical factorization
 $(\ref{SpectrRozclad_f_cont})$  is the least
favorable in the class of admissible spectral densities
$\mathcal{D}$ for the optimal linear estimation of the functional
$A\xi$ if
\begin{equation}\label{minmax f_cont}f^0(\lambda)=\left|\int_{0}^{\infty}\varphi^0(t)e^{-i\lambda t}dt\right|^2, \end{equation}
 where $\varphi^0(t)$, $t\in[0;\infty)$ is a solution to the conditional extremum problem
\begin{equation}||\mathbf D^{\tau}\mathbf A\varphi_{\tau}||^2\to\max,\quad
f(\lambda)=\left|\int_{0}^{\infty}\varphi(t)e^{-i\lambda
t}dt\right|^2\in\mathcal{D}.\label{zadumextpextr_A_cont}\end{equation}
\end{lem}

\begin{lem} Spectral density  $f^0(\lambda)\in\mathcal{D}$ which admits the
canonical factorization $(\ref{SpectrRozclad_f_cont})$  is the least
favorable in the class of admissible spectral densities $\mathcal{D}$
for the optimal linear estimation of the functional $A_T\xi$ if
 \begin{equation}\label{minmax f_T_cont}f^0(\lambda)=\left|\int_{0}^{T}\varphi^0(t)e^{-i\lambda t}dt\right|^2, \end{equation}
 where $\varphi^0(t)$, $t\in[0;T]$, is a solution to the conditional extremum problem
\begin{equation}||\mathbf D^{\tau}_T\mathbf A_T\varphi_{\tau}||^2\to\max,\quad
f(\lambda)=\left|\int_{0}^{T}\varphi(t)e^{-i\lambda
t}dt\right|^2\in\mathcal{D}.\label{zadumextpextr_A_T_cont}\end{equation}
\end{lem}

If $h_{\tau}^ {(a)}(f^0)\in H_{\mathcal{D}}$, the
minimax-robust spectral characteristic can be calculated as
$h^0=h_{\tau}^{(a)}(f^0)$.

The minimax-robust spectral characteristic $h^0$ and the least
favorable spectral density $f^0$ form a saddle point of the function $\Delta(h;f)$
on the set $H_{\mathcal{D}}\times\mathcal{D}$. The
saddle point inequalities
\[\Delta(h;f^0)\geq\Delta(h^0;f^0)\geq\Delta(h^0;f)\quad\forall f\in
\mathcal{D},\forall h\in H_{\mathcal{D}}\]
 hold true if
$h^0=h_{\tau}^{(a)}(f^0)$ and $h_{\tau}^{(a)}(f^0)\in
H_{\mathcal{D}}$, where $f^0$ is a solution
to the conditional extremum problem
\begin{equation}
\label{extrem1_cont}
 \widetilde{\Delta}(f)=
-\Delta(h_{\tau}^{(a)}(f^0);f) \to\inf, \quad{f\in \mathcal{D}},
\end{equation}
\[
\Delta(h_{\tau}^{(a)}(f^0);f)=\frac{1}{2\pi}\int_{-\infty}^{\infty}
\frac{|r_{\tau }( \lambda )|^2}
{f^0(\lambda)}f(\lambda)d\lambda.\]
Here $r_{\tau}$ is determined by formula $(\ref{r_A_cont})$ or $(\ref{r_A_T_cont})$ with $f(\lambda)=f^0(\lambda)$.
  The conditional extremum problem (\ref{extrem1_cont}) is equivalent to
the unconditional extremum problem
\[\Delta_{\mathcal{D}}(f)=\widetilde{\Delta}(f)+\delta(f|\mathcal{D})\to\inf,\]
where $\delta(f|\mathcal   D)$ is the indicator function
of the set $\mathcal D$. Solution $f^0$ to this unconditional
extremum problem is characterized by the condition $0\in
\partial\Delta_{\mathcal{D}}(f^0)$ (see Pshenichnyi\cite{Pshenychn}),
where $\partial\Delta_{\mathcal{D}}(f^0)$ is the subdifferential of
the functional $\Delta_{\mathcal{D}}(f^0)$ at point $f^0$. With the help of the
condition $0\in
\partial\Delta_{\mathcal{D}}(f^0)$ we can find the least
favorable spectral densities in some special classes of spectral
densities (see books by Moklyachuk\cite{Moklyachuk:2008}, Moklyachuk and Masyutka\cite{Moklyachuk:2012}
for more details).

\section{\Large{Least favorable spectral densities in the class $\mathcal{D}_0$}}
Consider the problem of the optimal estimation of functionals $A\xi$
and $A_T\xi$ of unknown values $\xi(t)$, $t\geq0$, of the
random process $\xi(t)$ with stationary $n$th increments in the
case where the spectral density is not known, but the following set
of spectral densities is given
$$\mathcal{D}_0=\left\{f(\lambda)|\frac{1}{2\pi}\int_{-\infty}^{\infty} f(\lambda)d\lambda\leq P_0\right\}.$$
It follows from the condition $0\in
\partial\Delta_{\mathcal{D}}(f^0)$ for $\mathcal{D}=\mathcal{D}_0$ that  the least favorable density satisfies the equation
$$|r_{\tau}^{(a)}( \lambda )|^2(f^0(\lambda))^{-1}=\psi(\lambda)+c^{-2},$$
where $\psi(\lambda)\leq0$ and $\psi(\lambda)=0$ if $f^0(\lambda)>0$. Therefore, the least favorable density in the
class $\mathcal{D}_0$ for the optimal linear estimation of the
functional $A\xi$ can be presented in the form
\begin{equation}f^0(\lambda)=\left|c\int_{0}^{\infty}\mathbf D^{\tau} (\mathbf A\varphi_{\tau}^0)(t)e^{i\lambda t}dt
\right|^2,\label{minmaxD0_cont}\end{equation} where the unknown function
$c\varphi_{\tau}^0(t)$ can be calculated using factorization
$(\ref{SpectrRozclad_f_cont})$, equation $(\ref{dd1_cont})$, condition
$(\ref{zadumextpextr_A_cont})$ and condition
$\int_{-\infty}^{\infty} |\varphi^0(\lambda)|^2d\lambda=2\pi P_0.$

Consider the equation
 \begin{equation}\mathbf D^{\tau}\mathbf A\mathbf W^{\tau}\varphi=\alpha \overline{\varphi},\quad
 \alpha\in  C.\label{rivnVl_cont}\end{equation}   For each solution of this equation such that $||\varphi ||^2=\dfrac{1}{2\pi}\int_{-\infty}^{\infty}|\varphi(\lambda)|^2d\lambda=P_0$ the following relation holds true:  \[f^0(\lambda)=\left|\int_{0}^{\infty}\varphi(t)e^{-i\lambda t}dt\right|^2=\]\[=
 \left|c\int_{0}^{\infty}\mathbf D^{\tau} (\mathbf A\mathbf W^{\tau}\varphi)(t)e^{i\lambda t}dt\right|^2.\]
Denote by $\nu_0P_0$ the maximum value of $||\mathbf D^{\tau}\mathbf A\mathbf W^{\tau}\varphi||^2$ on the set of those solutions
$\varphi$ of  equation $(\ref{rivnVl_cont})$, which satisfy condition
$||\varphi||^2=P_0$ and define canonical factorization
$(\ref{SpectrRozclad_f_cont})$ of the spectral density $f^0(\lambda)$.
Let $\nu_0^+P_0$ be the maximum value of $||\mathbf D^{\tau}\mathbf A\mathbf W^{\tau}\varphi||^2$ on the set of those $\varphi$ which
satisfy condition $||\varphi||^2=P_0$ and define  canonical
factorization $(\ref{SpectrRozclad_f_cont})$ of the spectral density
$f^0(\lambda)$ defined by $(\ref{minmaxD0_cont})$.

The derived equations and conditions give us a possibility to verify the validity of
following statement.

\begin{thm} If there exists a solution
$\varphi^0=\varphi^0(t)$ of equation $(\ref{rivnVl_cont})$
which satisfies conditions $||\varphi^0||^2=P_0$
and $\nu_0P_0=\nu_0^+P_0=||\mathbf D^{\tau}\mathbf A\mathbf W^{\tau}\varphi^0||^2$, the spectral
density $(\ref{minmax f_cont})$ is the least favorable density in the class
  $\mathcal{D}_0$ for the optimal estimation of the functional $A\xi$
of unknown values $\xi(t)$, $t\geq0$, of the random process
$\xi(t)$ with stationary
 $n$th increments. The increment $\xi^{(n)}(t,\tau)$ admits a one-sided moving average representation.
 If $\nu_0<\nu_0^+$, the density $(\ref{minmaxD0_cont})$ which admits the canonical factorization
$(\ref{SpectrRozclad_f_cont})$ is the least favorable in the class
$\mathcal{D}_0$. The function
 $c\varphi_{\tau}=c\varphi_{\tau}(t)$ is determined by equality $(\ref{dd1_cont})$, condition
$(\ref{zadumextpextr_A_cont})$ and condition $\int_{-\infty}^{\infty} |\varphi(\lambda)|^2d\lambda=2\pi
P_0$. The minimax-robust spectral
characteristic is calculated by formulas $(\ref{spectr_har_A_cont})$,
$(\ref{r_A_cont})$ substituting $f(\lambda)$ by $f^0(\lambda)$. \end{thm}

Consider the problem of optimal estimation of the functional $A_T\xi$. In this case the least favorable spectral density is
determined by the relation \begin{equation}f^0(\lambda)=\left|c\int_{0}^{T}\mathbf D^{\tau}_T (\mathbf A_T\varphi_{\tau}^0)(t)e^{i\lambda t}dt\right|^2.\label{minmaxD0_T_cont}\end{equation} Define a linear operator  $ \widehat{\mathbf A}_T$ in $L_2([0,\infty))$ space by relation  \begin{equation}\label{operator A_T_cont}(\widehat{\mathbf A}_T\varphi_{\tau})(t)=\int_0^ta(T-t+u)\varphi_{\tau}(u)du.\end{equation}
 Taking into consideration $(\ref{dd1_cont})$, we have the following equality \[\left|r_{\tau,T}^{(a)}(\lambda)\right|^2=
\left|\int_{0}^{T}\mathbf D^{\tau}_T (\mathbf A_T\mathbf W^{\tau}\varphi)(t)e^{i\lambda t}dt\right|^2=\]\begin{equation}=
\left|\int_{0}^{T}\mathbf D^{\tau}_T ( \widehat{\mathbf A}_T\mathbf W^{\tau}\varphi)(t)e^{-i\lambda t}dt\right|^2,\label{spivv N_cont}\end{equation} where the linear operator $\mathbf D^{\tau}_T \widehat{\mathbf A}_T\varphi_{\tau}$ in $L_2([0,\infty))$ space is culculated by formula \[\mathbf D^{\tau}_T (\widehat{\mathbf A}_T\varphi_{\tau})(t)=\sum_{k=0}^{\left[\frac{t}{\tau}\right]}\int_0^{t-\tau k}a(T-t+u+\tau k)\varphi_{\tau}(u)d(k)du.\]

Therefore each solution $\varphi =\varphi(t)$, $t\in[0,T]$, of
the equation \begin{equation}
\label{rivnVl_T1_cont}\mathbf D^{\tau}_T \mathbf A_T\mathbf W^{\tau}\varphi=\alpha \overline{\varphi},\quad
 \alpha\in  C,\end{equation}
 or the equation
 \begin{equation}\label{rivnVl_T2_cont}\mathbf D^{\tau}_T  \widehat{\mathbf A}_T\mathbf W^{\tau}\varphi=\beta \overline{\varphi},\quad
 \beta\in  C,\end{equation} such that $||\varphi||^2=P_0$, satisfies the following equality  \[f^0(\lambda)=\left|\int_{0}^{N}\varphi(t)e^{-i\lambda t}dt\right|^2=
 \left|cr_{\tau,T}^{(a)}(\lambda)\right|^2. \]

Denote by $\nu_0^TP_0$  the maximum value of $||\mathbf D^{\tau}_T\mathbf A_T\mathbf W^{\tau}\varphi||^2=||\mathbf D^{\tau}_T\widehat{\mathbf A}_T\mathbf W^{\tau}\varphi||^2$ on the set of solutions $\varphi_T$ of the  equation $(\ref{rivnVl_T1_cont})$
or the equation $(\ref{rivnVl_T2_cont})$, which satisfy condition
$||\varphi||^2=P_0$ and determine the  canonical factorization
$(\ref{SpectrRozclad_f_cont})$ of the spectral density $f^0(\lambda)\in
\mathcal{D}_0$. Let $\nu_0^{N+}P_0$ be the maximum value of
  $||\mathbf D^{\tau}_T\mathbf A_T\mathbf W^{\tau}\varphi||^2$ on the set of those
$\varphi$ which satisfy condition $||\varphi||^2=P_0$ and
determine the  canonical factorization ($\ref{SpectrRozclad_f_cont}$) of
the spectral density $f^0(\lambda)$ defined by ($\ref{minmaxD0_T_cont}$).

The following statement holds true.

\begin{thm} If there exists a solution
$\varphi^0=\varphi^0(t)$, $t\in[0;T]$ of  equation $(\ref{rivnVl_T1_cont})$ or  equation $(\ref{rivnVl_T2_cont})$ such that $||\varphi^0||^2=P_0$
and $\nu_0^TP_0=\nu_0^{ T+}P_0=||\mathbf D^{\tau}_T\mathbf A_T\mathbf W^{\tau}\varphi^0||^2$, the spectral density $(\ref{minmax f_T_cont})$ is  least favorable in the
class $\mathcal{D}_0$ for the optimal estimation of the functional
$A_T\xi$ of unknown values $\xi(t)$, $t\in[0;T]$, of the
random process $\xi(t)$ with stationary
 $n$th increments. The increment $\xi^{(n)}(t,\tau)$ admits a one-sided moving average representation. If $\nu_0^T<\nu^{T+}_0 $, the density $(\ref{minmaxD0_T_cont})$
which admits the canonical factorization $(\ref{SpectrRozclad_f_cont})$
is the least favorable in the class $\mathcal{D}_0$. The  function
 $c\varphi_{\tau}=c\varphi_{\tau }(t)$, $t\in[0;T]$, is determined by equation $(\ref{dd1_cont})$, condition
$(\ref{zadumextpextr_A_T_cont})$ and condition $\int_{-\infty}^{\infty} |\varphi(\lambda)|^2d\lambda=2\pi
P_0$. The minimax-robust spectral characteristic is calculated by formulas
$(\ref{spectr_har_A_T_cont})$, $(\ref{r_A_T_cont})$ substituting $f(\lambda)$ by $f^0(\lambda)$.\end{thm}

\section{\Large{Least favorable spectral densities in the class $\mathcal{D}_v^u$}}

Consider the case where the spectral density $f(\lambda)$ is not known, but the following set
of spectral densities is given: \[\mathcal{D}^u_v=\left \{f(\lambda)|v(\lambda)\leq f(\lambda)\leq u(\lambda),\,\int_{-\infty}^{\infty} f(\lambda)d\lambda=2\pi P_0\right \}\] where $v(\lambda)$ and $u(\lambda)$ are some given (fixed)
spectral densities. It follows from the condition $0\in
\partial\Delta_{\mathcal{D}}(f^0)$ for $\mathcal{D}=\mathcal{D}^u_v$ that the least favorable density $f^0(\lambda)$ in the class $\mathcal{D}_v^u$ for the optimal linear estimation of the functional
$A\xi$ is of the form
\begin{equation}f^0(\lambda)=\max\left \{v(\lambda),\min\left \{u(\lambda),\left |s^0_{\tau}(\lambda)
\right
|^2\right \}\right \},\label{minmaxDuv_cont}\end{equation}
\[s^0_{\tau}(\lambda)=c\int_{0}^{\infty}\mathbf D^{\tau} (\mathbf A\varphi_{\tau}^0)(t)e^{i\lambda t}dt,\]
 where the unknown function
$c\varphi_{\tau}^0(t)$ can be calculated using factorization
$(\ref{SpectrRozclad_f_cont})$, equestion $(\ref{dd1_cont})$, conditions
$(\ref{zadumextpextr_A_cont})$ and
$\int_{-\infty}^{\infty} |\varphi^0(\lambda)|^2d\lambda=2\pi P_0.$

Denote by $\nu_{uv}P_0$ the maximum value of $||\mathbf D^{\tau}\mathbf A\mathbf W^{\tau}\varphi||^2$ on the set of those solutions
$\varphi$ of  equation $(\ref{rivnVl_cont})$ which satisfy condition
$||\varphi||^2=P_0$, inequalities
\[v(\lambda)\leq\left |\int_{0}^{\infty}\varphi(t)e^{-i\lambda t}dt\right |^2\leq u(\lambda)\]
and determine the canonical
factorization $(\ref{SpectrRozclad_f_cont})$ of the spectral density $f(\lambda)$. Let  $\nu_{uv}^+P_0$ be the maximum value of
  $||\mathbf D^{\tau}\mathbf A\mathbf W^{\tau}\varphi||^2$ on the set of those $\varphi$ which
satisfy condition $||\varphi||^2=P_0$ and determine the canonical
factorization $(\ref{SpectrRozclad_f_cont})$ of the spectral density
$f^0(\lambda)$ defined by $(\ref{minmaxDuv_cont})$.

 The derived equations and conditions give us a possibility to verify
the validity of the following statement.

\begin{thm} If there exists a solution
$\varphi^0=\varphi^0(t)$ of equation $(\ref{rivnVl_cont})$
which satisfies conditions $||\varphi^0||^2=P_0$
and $\nu_{uv}P_0=\nu_{uv}^+P_0=||\mathbf D^{\tau}\mathbf A\mathbf W^{\tau}\varphi^0||^2$, the spectral
density $(\ref{minmax f_cont})$ is the least favorable in the class
  $\mathcal{D}_v^u$ for the optimal estimation of the functional $A\xi$
of unknown values $\xi(t)$, $t\geq0$, of the random process
$\xi(t)$ with stationary
 $n$th increments.
 The increment $\xi^{(n)}(t,\tau)$ admits one-sided moving average representation. If $\nu_{uv}<\nu_{uv}^+$, the density $(\ref{minmaxDuv_cont})$ which admits the canonical factorization
$(\ref{SpectrRozclad_f_cont})$ is the least favorable in the class
$\mathcal{D}_v^u$. The function
 $c\varphi_{\tau}=c\varphi_{\tau}(t)$ is determined by equality $(\ref{dd1_cont})$, conditions
$(\ref{zadumextpextr_A_cont})$ and $\int_{-\infty}^{\infty} |\varphi(\lambda)|^2d\lambda=2\pi
P_0$. The minimax-robust spectral
characteristic is calculated by formulas $(\ref{spectr_har_A_cont})$,
$(\ref{r_A_cont})$ substituting $f(\lambda)$ by $f^0(\lambda)$.\end{thm}

Consider the problem of the optimal estimation of the functional $A_T\xi$. In this case the least favorable spectral density is
determined by the relation
\begin{equation}f^0(\lambda)=\max\left \{v(\lambda),\min\left \{u(\lambda),\left |s^0_{\tau,T}(\lambda)
\right |^2\right \}\right \}.\label{minmaxDuv_T_cont}
\end{equation}
\[s^0_{\tau,T}(\lambda)=c\int_{0}^{T}\mathbf D^{\tau}_T (\mathbf A_T\varphi_{\tau})(t)e^{i\lambda t}dt\]
Denote  by $\nu_{uv}^TP_0$ the maximum value of $||\mathbf D^{\tau}_T\mathbf A_T\mathbf W^{\tau}\varphi||^2=||\mathbf D^{\tau}_T\widehat{\mathbf A}_T\mathbf W^{\tau}\varphi||^2$ on the set of those solutions $\varphi$ of equations  $(\ref{rivnVl_T1_cont})$ and $(\ref{rivnVl_T2_cont})$, which satisfy condition
$||\varphi||^2=P_0$,  inequalities
\[v(\lambda)\leq\left |\int_{0}^{T}\varphi(t)e^{-i\lambda t}dt\right |^2\leq u(\lambda)\]  and define the
canonical factorization $(\ref{SpectrRozclad_f_cont})$ of the spectral
density $f(\lambda)$. Let $\nu_{uv}^{T+}P_0$ be the
maximum value of $||\mathbf D^{\tau}_T\mathbf A_T\mathbf W^{\tau}\varphi||^2$ on the set
of those $\varphi$ which satisfy condition
$||\varphi||^2=P_0$ and define canonical factorization $(\ref{SpectrRozclad_f_cont})$ of the
spectral density $f^0(\lambda)$ determined by  $(\ref{minmaxDuv_T_cont})$.

The following statement holds true.

\begin{thm} If there exists a solution
$\varphi^0=\varphi^0(t)$, $t\in[0;T]$, of equation $(\ref{rivnVl_T1_cont})$ or equation $(\ref{rivnVl_T2_cont})$ which satisfies conditions $||\varphi^0||^2=P_0$
and $\nu_{uv}^TP_0=\nu_{uv}^{ T+}P_0=||\mathbf D^{\tau}_T\mathbf A_T\mathbf W^{\tau}\varphi^0||^2$, spectral density
$(\ref{minmax f_T_cont})$ is  least favorable in
the class $\mathcal{D}_v^u$ for the optimal estimation of the
functional $A_T\xi$
of unknown values $\xi(t)$, $t\in[0;T]$, of the random process $\xi(t)$ with stationary
 $n$th increments. The increment $\xi^{(n)}(t,\tau)$ admits one-sided moving average representation. If $\nu_{uv}^T<\nu_{uv}^{T+}$, the density $(\ref{minmaxDuv_T_cont})$,
which admits the canonical factorization $(\ref{SpectrRozclad_f_cont})$,
is  least favorable in the class $\mathcal{D}_v^u$. The function
 $c\varphi_{\tau}=c\varphi_{\tau }(t)$, $t\in[0;T]$, is determined by equation $(\ref{dd1_cont})$, conditions
$(\ref{zadumextpextr_A_T_cont})$ and $\int_{-\infty}^{\infty} |\varphi(\lambda)|^2d\lambda=2\pi
P_0$. The minimax-robust spectral characteristic is calculated by formulas
$(\ref{spectr_har_A_T_cont})$, $(\ref{r_A_T_cont})$ substituting $f(\lambda)$ by $f^0(\lambda)$.\end{thm}

\section{\Large{Least favorable spectral densities in the class $\mathcal{D}_{\delta}$}}

Consider the problem of the optimal estimation of the functionals $A\xi$
and $A_T\xi$ of unknown values $\xi(t)$, $t\geq0$, of the
random process $\xi(t)$ with stationary $n$th increments in the
case where the spectral density is not known, but the following set
of spectral densities is given
 \[\mathcal{D}_{\delta}=\left\{f(\lambda)| \frac{1}{2\pi}\int_{-\infty}^{\infty} |f(\lambda)-v(\lambda)|d\lambda\leq \delta\right\},\] where $v(\lambda)$ is a bounded spectral density.
It comes from the condition $0\in
\partial\Delta_{\mathcal{D}_{\delta}}(f^0)$ that    the least favorable
spectral densities in the class $\mathcal{D}_{\delta}$ for optimal linear extrapolation of the functional $A\xi$
can be presented in the form
\begin{equation}f^0(\lambda)=\max\left \{v(\lambda),\left |
c\int_{0}^{\infty}\mathbf D^{\tau} (\mathbf A\varphi_{\tau}^0)(t)e^{i\lambda t}dt\right
|^2\right \},\label{minmaxDdelta_cont}\end{equation} where unknown function
$c\varphi_{\tau}^0(t)$ is calculated using the factorization
$(\ref{SpectrRozclad_f_cont})$, relation $(\ref{dd1_cont})$, condition
$(\ref{zadumextpextr_A_cont})$ and condition
\[\frac{1}{2\pi}\int_{-\infty}^{\infty} |\varphi^0(\lambda)|^2d\lambda=\delta+\frac{1}{2\pi}\int_{-\infty}^{\infty} v(\lambda)d\lambda=P_1.\]

Define by $\nu_{\delta}P_1$ the maximum value of $||\mathbf D^{\tau}\mathbf A\mathbf W^{\tau}\varphi||^2$ on the set of those $\varphi$ which
belongs to the set of solutions of equation
  $(\ref{rivnVl_cont})$, satisfy equation
$||\varphi||^2=P_1$, inequality
\[v(\lambda)\leq\left |\int_{0}^{\infty}\varphi(t)e^{-i\lambda t}dt\right |^2\]
and determine the canonical
factorization $(\ref{SpectrRozclad_f_cont})$
of the spectral density $f(\lambda)$. Let $\nu_{\delta}^+P_1$ be the maximum value of $||\mathbf D^{\tau}\mathbf A\mathbf W^{\tau}\varphi||^2$ on the set of those $\varphi$, which
satisfy condition $||\varphi||^2=P_1$ and determine the canonical
factorization $(\ref{SpectrRozclad_f_cont})$ of the spectral density
$f^0(\lambda)$ defined by $(\ref{minmaxDdelta_cont})$. The following
statement holds true.

\begin{thm} If there exists a solution
$\varphi^0=\varphi^0(t)$ of equation $(\ref{rivnVl_cont})$
which satisfies conditions $||\varphi^0||^2=P_1$
and $\nu_{\delta}P_0=\nu_{\delta}^+P_1=||\mathbf D^{\tau}\mathbf A\mathbf W^{\tau}\varphi^0||^2$, the spectral density $(\ref{minmax f_cont})$ is the least favorable in the
class
$\mathcal{D}_{\delta}$ for the optimal extrapolation of the functional $A\xi$
of unknown values $\xi(t)$, $t\geq0$, of the random process
$\xi(t)$ with stationary
 $n$th increments. The increment $\xi^{(n)}(t,\tau)$ admits one-sided moving average representation. If  $\nu_{\delta}<\nu_{\delta}^+$, the density $(\ref{minmaxDdelta_cont})$,
 which admits the canonical factorization $(\ref{SpectrRozclad_f_cont})$,
is  least favorable  in the class $\mathcal{D}_{\delta}$. The function
 $c\varphi_{\tau}=c\varphi_{\tau}(t)$ is determined by equality $(\ref{dd1_cont})$, condition
$(\ref{zadumextpextr_A_cont})$ and condition
\begin{equation}\frac{1}{2\pi}\int_{-\infty}^{\infty} |\varphi(\lambda)|^2d\lambda=\delta+\frac{1}{2\pi}\int_{-\infty}^{\infty} v(\lambda)d\lambda=P_1.\label{definition_p2_cont}
\end{equation}
The minimax-robust spectral
characteristic is calculated by formulas $(\ref{spectr_har_A_cont})$,
$(\ref{r_A_cont})$ substituting $f(\lambda)$ by $f^0(\lambda)$.\end{thm}

In the case of  optimal estimation of the functional $A_T\xi$ the
least favorable spectral density is determined by formula
\begin{equation}f^0(\lambda)=\max\left \{v(\lambda),\left |
c\int_{0}^{T}\mathbf D^{\tau}_T (\mathbf A_T\varphi_{\tau}^0)(t)e^{i\lambda t}dt\right |^2\right \}.\label{minmaxDdelta_T_cont}
\end{equation}

Let $\nu_{\delta}^TP_1$ be the maximum value of $||\mathbf D^{\tau}_T\mathbf A_T\mathbf W^{\tau}\varphi||^2=||\mathbf D^{\tau}_T\widehat{\mathbf A}_T\mathbf W^{\tau}\varphi||^2$ on the set of those $\varphi$ which belong to the set of solutions
of equation $(\ref{rivnVl_T1_cont})$ and $(\ref{rivnVl_T2_cont})$, satisfy condition
$||\varphi||^2=P_1$, the inequality
\[v(\lambda)\leq\left |\int_{0}^{T}\varphi(t)e^{-i\lambda t}dt\right |^2\]  and determined the
canonical factorization $(\ref{SpectrRozclad_f_cont})$ of the spectral
density $f(\lambda)$. Let $\nu_{\delta}^{T+}P_1$ be the maximum value of $||\mathbf D^{\tau}_T\mathbf A_T\mathbf W^{\tau}\varphi||^2$ on the set of those  $\varphi$ which
satisfy condition
$||\varphi||^2=P_1$ and determined the canonical
factorization $(\ref{SpectrRozclad_f_cont})$ of the spectral density
$f^0(\lambda)$ defined by $(\ref{minmaxDdelta_T_cont})$. The following
statement holds true.

\begin{thm} If there exists a solution
$\varphi^0=\varphi^0(t)$, $t\in[0;T]$, of equation $(\ref{rivnVl_T1_cont})$ or equation $(\ref{rivnVl_T2_cont})$
 which satisfies conditions $||\varphi^0||^2=P_1$
and $\nu_{\delta}^TP_1=\nu_{\delta}^{ T+}P_1=||\mathbf D^{\tau}_T\mathbf A_T\mathbf W^{\tau}\varphi^0||^2$, the spectral density $(\ref{minmax f_T_cont})$ is  least favorable  in
the class
$\mathcal{D}_{\delta}$ for the optimal extrapolation
of the functional $A_T\xi$ of unknown values $\xi(t)$, $t\in[0;T]$, of the random process $\xi(t)$ with
stationary
 $n$th increments. The increment $\xi^{(n)}(t,\tau)$ admits one-sided moving average representation. If $\nu_{\delta}^T<\nu_{\delta}^{T+}$,  the density $(\ref{minmaxDdelta_T_cont})$
which admits the canonical factorization $(\ref{SpectrRozclad_f_cont})$
is the least favorable in the class  $\mathcal{D}_{\delta}$. The function
 $c\varphi_{\tau}=c\varphi_{\tau }(t)$, $t\in[0;T]$, is determined by equation $(\ref{dd1_cont})$, conditions
$(\ref{zadumextpextr_A_T_cont})$ and $(\ref{definition_p2_cont})$. The minimax-robust spectral
characteristic is calculated by formulas $(\ref{spectr_har_A_T_cont})$,
$(\ref{r_A_T_cont})$  substituting $f(\lambda)$ by $f^0(\lambda)$.\end{thm}

\section{\Large{Conclusions}}

In this article methods of solution of the problem of
optimal linear estimation of functionals
\[A {\xi}=\int_{0}^{\infty}a(t)\xi(t)dt,\quad
A_T{\xi}=\int_{0}^Ta(t)\xi(t)dt\]
which depend on unknown
values of a random process $\xi(t)$ with stationary $n$th
increments were described. The received estimates are based on observations of the process
$\xi(t)$ for $t<0$. Formulas are derived for
computing the value of the mean-square error and the spectral
characteristic of the optimal linear estimates of functionals in the
case of spectral certainty where the spectral density of the
process is known.

In the case of spectral uncertainty where the spectral density is
not known but, instead, a set of admissible spectral densities is
specified, the minimax-robust method was applied. We propose a
representation of the mean square error in the form of a linear
functional in $L_1$ with respect to spectral density, which allows
us to solve the corresponding conditional extremum problem and
describe the minimax (robust) estimates of the functionals. Formulas
that determine the least favorable spectral densities and minimax
(robust) spectral characteristics of the optimal linear estimates of
the functionals are derived in this case for some concrete classes
of admissible spectral densities.

\noindent\hrulefill

\end{document}